\documentstyle{article}
\begin{document}
\baselineskip+8pt \small \begin{center}{\textit{ In the name of
 Allah, the Beneficent, the Merciful}}\end{center}

\begin{center}{\bf  A matrix representation of composition of polynomial maps}
\end{center}
\begin{center}{\bf  Ural Bekbaev  \\
 INSPEM, Universiti Putra Malaysia.\\
e-mail: bekbaev@science.upm.edu.my}
 \end{center}

 \begin{abstract}{In this paper polynomial maps are represented by the use of matrices whose entries are numbered by pair of
 multiindices.  A new product of such matrices is introduced. By the use of this and ordinary product of
 matrices the matrix representation of composition of polynomial maps is
 given. In the case of real and complex numbers different kind of norms of such matrices are introduced.
 Properties of these norms with respect to the ordinary and new products are investigated. A generalization of
 Bombieri's inequality is offered.}\end{abstract}


  \vspace{0.5cm}

 In this paper we are going to offer a matrix representation for polynomial maps and their compositions.
 For this purpose a new product of matrices, whose  entries are numbered by pair of
 multiindices, is introduced. The matrix representation of composition of two polynomial maps
 is given. In the case of real and complex numbers different kind of norms of such matrices are introduced.
 Properties of these norms with respect to the ordinary and new products are investigated. A generalization of
 Bombieri's inequality is offered.

 For a positive integer $n$ let $I_n$ stand for all row $n$-tuples with nonnegative integer entries with the
 following linear order: $\beta=
 (\beta_1,\beta_2,...,\beta_n)<\alpha=(\alpha_1,\alpha_2,...,\alpha_n)$ if and only if $\vert \beta\vert < \vert \alpha\vert$ or
$\vert \beta\vert = \vert \alpha\vert$ and $\beta_1> \alpha_1$ or
$\vert \beta\vert = \vert \alpha\vert$, $\beta_1= \alpha_1$ and
$\beta_2> \alpha_2$ etcetera, where $\vert \alpha\vert$ stands for
$\alpha_1+\alpha_2+...+\alpha_n $.

It is clear that for $\alpha, \beta, \gamma \in I_n$ one has
$\alpha < \beta $ if and only if $\alpha + \gamma < \beta
+\gamma$. We write $\beta \ll \alpha$ if $\beta_i \leq \alpha_i$
for all $i=1,2,...,n$, $\left(\begin{array}{c}
  \alpha \\
  \beta \\
\end{array}\right)$ stands for $\frac{\alpha!}{\beta!
(\alpha -\beta)!}$, $\alpha !=\alpha_1!\alpha_2!...\alpha_n!$.

In the case of $n=0$ it is assumed that $I_n=  \{0\}$.

In future we use the following result.

{\bf Proposition 1.}  If $\vert \alpha\vert =p+q$, where $p$, $q$
are nonnegative integer numbers,then
$$\sum_{\beta \ll \alpha, \vert \beta\vert =p}\left(
\begin{array}{c}
  \alpha \\
  \beta \\
\end{array}\right)= \left(
\begin{array}{c}
  p+q \\
  p \\
\end{array}\right)$$

 In future $n$, $n'$ and $n''$ are assumed to be any fixed nonnegative integers.

 For any commutative, associative ring $R$ containing the field of rational numbers and nonnegative
integer numbers $p,p'$ let  $M_{n,n'}(p,p';R)=M(p,p';R)$ stand for
all $"p\times p'"$ size matrices $A=(A_{\alpha,\alpha'})_{\vert
\alpha \vert=p, \vert \alpha'\vert=p'}$ ($\alpha$ presents
row,$\alpha'$ presents column and $\alpha\in I_n,\alpha'\in
I_{n'}$) with entries from $R$. Over such kind matrices in
addition to ordinary sum and product of matrices we consider the
following "product" as well:

{\bf Definition 1.} If $A\in M(p,p';R)$ and $B\in M(q,q';R)$ then
$A\bigodot B=C\in M(p+q,p'+q';R)$ that for any
$\vert\alpha\vert=p+q$, $\vert\alpha'\vert=p'+q'$, where
$\alpha\in I_n,\alpha'\in I_{n'}$,
$$C_{\alpha,\alpha'}=\sum_{\beta,\beta'}\left(\begin{array}{c}
  \alpha \\
  \beta \\
\end{array}\right)
    A_{\beta,\beta'}B_{\alpha-\beta,\alpha'-\beta'}
$$, where the sum is taken over all $\beta\in I_n,\beta'\in I_{n'}$, for which $\vert
\beta\vert=p$, $\vert \beta'\vert=p'$, $\beta\ll \alpha$ and
$\beta'\ll \alpha'$.

Let us agree that $h$ ($H$, $v$, $V$) stands for any element of
$M(0,1;R)$ (respect. $M(0,p;R)$, $M(1,0;R)$, $M(p,0;R)$ , where
$p$ may be any nonnegative integer). We use $E_k$ for $"k\times
k"$ size ordinary unit matrix from $M_{n,n}(k,k;R)$. For the sake
of convenience it will be assumed that $A_{\alpha,\alpha'}=0$
($\alpha !=\infty$) whenever $\alpha \notin I_n$ or $\alpha'
\notin I_{n'}$ (respect. $\alpha \notin I_n$).

{\bf Proposition 2.} For the above defined product the following
are true.

1. $A\bigodot B=B\bigodot A$.

2. $(A+B)\bigodot C=A\bigodot C+ B\bigodot C$.

3. $(A\bigodot B)\bigodot C=A\bigodot (B\bigodot C)$

4. $ (\lambda A)\bigodot B=\lambda (A\bigodot B)$ for any
$\lambda\in R$

5. If $R$ is an integral domain then $A\bigodot B=0$ if and only
    if $A=0$ or $B=0$.

6. $A(B\bigodot H)=(AB)\bigodot H$

7. $(E_k\bigodot V)A=A\bigodot V$

{\bf Proof}. Let us show only that if $R$ is an integral domain
then $C=A\bigodot B=0$ if and only if $A=0$ or $B=0$. Assume that
$A\neq 0$, $B\neq 0$ and let $\beta_0$ $(\gamma_0)$ be the least
multiindex for which the corresponding row of $A$ (respec. $B$) is
nonzero. Similarly, let $\beta'_0$ $(\gamma'_0)$ be the least
multiindex for which $A_{\beta_0,\beta'_0}$ (respec.
$B_{\gamma_0,\gamma'_0}$) is not zero. Consider $\alpha
=\beta_0+\gamma_0$, $\alpha' =\beta'_0+\gamma'_0$ and
$C_{\alpha,\alpha'}$. It is easy to see that
$C_{\alpha,\alpha'}=\left(\begin{array}{c}
  \alpha \\
  \beta_0 \\
\end{array}\right)A_{\beta_0,\beta'_0}B_{\gamma_0,\gamma'_0} \neq 0$.
 It is the needed contradiction to complete the proof.

Proofs of the following two propositions are not difficult.

{\bf Proposition 3.} If $A_i\in M(p_i,q_i;R)$ for $i=1,2,...,m$,
$\vert\alpha\vert= p_1+p_2+...+p_m$, $\vert\alpha'\vert=
p_1'+p_2'+...+p_m'$ then

$(A_1\bigodot A_2\bigodot...\bigodot A_m)_{\alpha,\alpha'}=\sum
\frac{\alpha !}{\beta!
\gamma!...\delta!}A_{1\beta,\beta'}A_{2\gamma,\gamma'}...
A_{m\delta,\delta'}$,

 where the sum is taken over all $\beta,
\gamma,...,\delta\in I_n; \beta', \gamma',...,\delta'\in I_{n'}$
for which $\beta+ \gamma +...+\delta=\alpha$, $\beta'+ \gamma'
+...+\delta'=\alpha'.$

In future $A^{(m)}$ means the $m$-th power of matrix $A$ with
respect to the new product.

{\bf Proposition 4.} If $h=(h_1,h_2,...,h_{n})\in M(0,1;R)$,
$v=(v_1,v_2,...,v_n)\in M(1,0;R)$, then
$$(h^{(m)})_{0,\alpha'}=\left(\begin{array}{c}
  m \\
  \alpha'\\
\end{array}\right)h^{\alpha'}, \hspace{1cm} (v^{(m)})_{\alpha,0}=m!v^{\alpha}$$, where $h^{\alpha}$ stands for
$h_1^{\alpha_1}h_2^{\alpha_2}...h_n^{\alpha_n}$

{\bf Proposition 5.} For any nonnegative integers $p$, $q$, $p'$,
$q'$ and matrices $A\in M_{n,n'}(p,p';R)$, $B\in
M_{n,n'}(q,q';R)$, $h=(h_1,h_2,...,h_n)\in M_{n,n}(0,1;R)$,
$v=(v_1,v_2,...,v_{n'})\in M_{n',n'}(1,0;R)$ the following
equalities
$$(\frac{h^{(p)}}{p!}A)\bigodot
(\frac{h^{(q)}}{q!}B)=\frac{h^{(p+q)}}{(p+q)!}(A\bigodot B),
\hspace{1cm} (A\frac{v^{(p')}}{p'!})\bigodot
(B\frac{v^{(q')}}{q'!})=(A\bigodot
B)\frac{v^{(p'+q')}}{(p'+q')!}$$
 are true.

 {\bf Proof.}
$$((\frac{h^{(p)}}{p!}A)\bigodot (\frac{h^{(q)}}{q!}B))_{0,\alpha'}=
    \sum_{\beta'}(\frac{h^{(p)}}{p!}A)_{0,\beta'}(\frac{h^{(q)}}{q!}B)_{0,\alpha'-\beta'}=$$
    $$\sum_{\beta'}\sum_{\xi}(\frac{h^{(p)}}{p!})_{0,\xi}A_{\xi,\beta'}
    \sum_{\eta}(\frac{h^{(q)}}{q!})_{0,\eta}B_{\eta,\alpha'-\beta'}=\sum_{\beta'}
    \sum_{\xi}\frac{h^{\xi}}{\xi!}A_{\xi,\beta'}
    \sum_{\eta}\frac{h^{\eta}}{\eta!}B_{\eta,\alpha'-\beta'}=$$
$$\sum_{\beta'}\sum_{\xi,\eta}\left(\begin{array}{c}
  \xi+\eta \\
  \xi \\
\end{array}\right)\frac{h^{\xi+\eta}}{(\xi+\eta)!}A_{\xi,\beta'}
    B_{\eta,\alpha'-\beta'}=\sum_{\vert \xi \vert=p+q}(\frac{h^{(p+q)}}{(p+q)!})_{0,\xi}
    \sum_{\vert \eta\vert=p,\beta'}\left(\begin{array}{c}
  \xi \\
  \eta \\
\end{array}\right)A_{\eta,\beta'}B_{\xi-\eta,\alpha'-\beta'}=$$
$$\sum_{\vert \xi \vert=p+q}(\frac{h^{(p+q)}}{(p+q)!})_{0,\xi}(A\bigodot B)_{\xi,\alpha'}=
(\frac{h^{(p+q)}}{(p+q)!}A\bigodot B)_{0,\alpha'}$$ The proof of
the second identity is similar.

{\bf Remark 1.} Due to the "duality" of two equalities in
Proposition 5 in future  we will consider only results dealing
with the first equality. Analogies of the presented results
dealing with the second equality can be obtained in a similar way.

From Proposition 5 the following more general result can be
derived.

{\bf Proposition 6.} For any nonnegative integers $p$, $q$, $p'$,
$q'$, $k$ and matrices $A\in M_{n,n'}(k,1;R)$, $B\in
M_{n',n''}(p,p';R)$, $C\in M_{n',n''}(q,q';R)$, the following
equality
$$(\frac{A^{(p)}}{p!}B)\bigodot
(\frac{A^{(q)}}{q!}C)=\frac{A^{(p+q)}}{(p+q)!}(B\bigodot C)$$ is
true.

 {\bf Proof.} Due to Proposition 5 for $h'\in Mat_{n,n'}(0,1)$ one has the equality
 $$(\frac{h'^{(p)}}{p!}B)\bigodot (\frac{h'^{(q)}}{q!}C)=\frac{h'^{(p+q)}}{(p+q)!}(B\bigodot
 C)$$
Substitution $\frac{h^{(k)}}{k!}A$ for $h'$, where $h\in
M_{n,n}(0,1;R)$ into this equality implies that
$$(\frac{(\frac{h^{(k)}}{k!}A)^{(p)}}{p!}B)\bigodot
(\frac{(\frac{h^{(k)}}{k!}A)^{(q)}}{q!}C)=\frac{(\frac{h^{(k)}}{k!}A)^{(p+q)}}{(p+q)!}(B\bigodot
C)$$

The left side of this equality equals to
$$(\frac{h^{(kp)}}{(kp)!}\frac{A^{(p)}B}{p!})\bigodot
(\frac{h^{(kq)}}{(kq)!}\frac{A^{(q)}C}{q!})=\frac{h^{(k(p+q))}}{(k(p+q))!}(\frac{A^{(p)}B}{p!}\bigodot\frac{A^{(q)}C}{q!})$$
, the right side equals to
$$\frac{h^{(k(p+q))}}{(k(p+q))!}\frac{A^{(p+q)}}{(p+q)!}(B\bigodot
C)$$ Therefore the conclusion of Proposition 6 is true.

 In future let $F$ stand for the field of real or complex numbers, $\rho \geq 1$
 be any fixed real number and $\varrho$ stand for the real number for which $\frac{1}{\rho}+\frac{1}{\varrho}=1$.
We consider the following $\rho$-norm of elements $A\in
M(p,p';F)$:

{\bf Definition 2.}
$$\|A\|=\|A\|_{\rho}=(\sum_{\alpha,\alpha'}\frac{\vert A_{\alpha,\alpha'}\vert ^{\rho}}{\alpha!(p!p'!)^{\rho
-1}})^{1/\rho}$$

{\bf Theorem 1.} 1. If $A,B\in Mat(p,p';F)$ and $\lambda \in
F$ then

a)$\|A\|=0$ if and only if $A=0$,

b)$\|\lambda A\| =|\lambda |\|A\|$,

c) $\|A+B\|\leq \|A\|+ \|B\|$.

2. For any nonnegative integer numbers $p$, $p'$, $q$ and $q'$ there is such positive number $\lambda(p, p', q, q')$
that for any $A\in Mat(p,p';F)$, $B\in Mat(q,q';F)$ the
following inequality is valid:

 $$\lambda(p, p', q, q')\|A\|\|B\| \leq \|A\bigodot B\|\leq \|A\|\|B\|$$

{\bf Proof.}  Here is a proof of part 2. First let us show the
inequality $\|A\bigodot B\|\leq \|A\|\|B\|$.

     Due to the H\"{o}lder inequality  for $A\bigodot B=C$ one has
$$\vert C_{\alpha,\alpha'}\vert =\vert\sum_{\beta <<\alpha,\beta'<<\alpha'}\left(\begin{array}{c}
  \alpha \\
  \beta \\
\end{array}\right)A_{\beta,\beta'}B_{\alpha-\beta,\alpha'-\beta'}\vert \leq
\alpha !\sum_{\beta <<\alpha,\beta'<<\alpha'}\frac{\vert
A_{\beta,\beta'}B_{\alpha-\beta,\alpha'-\beta'}\vert}{(\beta
!(\alpha-\beta)!)^{1/\rho}}\frac{1}{(\beta
!(\alpha-\beta)!)^{1-1/\rho}} \leq $$ $$\alpha !(\sum_{\beta
<<\alpha,\beta'<<\alpha'}\frac{\vert
A_{\beta,\beta'}B_{\alpha-\beta,\alpha'-\beta'}\vert^{\rho}}{\beta
!(\alpha-\beta)!})^{1/\rho}(\sum_{\beta
<<\alpha,\beta'<<\alpha'}(\frac{1}{(\beta
!(\alpha-\beta)!)^{1-1/\rho}})^{\varrho})^{1/\varrho}=$$
$$\alpha !(\sum_{\beta
<<\alpha,\beta'<<\alpha'}\frac{\vert
A_{\beta,\beta'}\vert^{\rho}}{\beta !}\frac{\vert
B_{\alpha-\beta,\alpha'-\beta'}\vert^{\rho}}{(\alpha-\beta)!})^{1/\rho}(\sum_{\beta
<<\alpha}\frac{1}{\beta
!(\alpha-\beta)!}\sum_{\beta'<<\alpha'}1)^{1/\varrho} \leq$$
$$\alpha !(\sum_{\beta
<<\alpha,\beta'<<\alpha'}\frac{\vert
A_{\beta,\beta'}\vert^{\rho}}{\beta !}\frac{\vert
B_{\alpha-\beta,\alpha'-\beta'}\vert^{\rho}}{(\alpha-\beta)!})^{1/\rho}(\left(
\begin{array}{c}
  p+q \\
  p \\
\end{array}\right)\frac{1}{\alpha !}\left(
\begin{array}{c}
  p'+q' \\
  p' \\
\end{array}\right))^{1/\varrho}$$ as far as according to
Proposition 1 one has $\sum_{\beta <<\alpha}\frac{1}{\beta
!(\alpha-\beta)!}=\left(
\begin{array}{c}
  p+q \\
  p \\
\end{array}\right)\frac{1}{\alpha !}$ and
$\sum_{\beta'<<\alpha'}1\leq \left(
\begin{array}{c}
  p'+q' \\
  p' \\
\end{array}\right)$. Therefore
$$ \|C\|=(\sum_{\alpha,\alpha'}\frac{\vert C_{\alpha,\alpha'}\vert^{\rho}}{\alpha
!((p+q)!(p'+q')!)^{\rho-1}})^{1/\rho} \leq $$
$$(\sum_{\alpha,\alpha'}\frac{1}{{\alpha
!((p+q)!(p'+q')!)^{\rho-1}}}(\alpha !)^{\rho}\sum_{\beta
<<\alpha,\beta'<<\alpha'}\frac{\vert
A_{\beta,\beta'}\vert^{\rho}}{\beta !}\frac{\vert
B_{\alpha-\beta,\alpha'-\beta'}\vert^{\rho}}{(\alpha-\beta)!}(\left(
\begin{array}{c}
  p+q \\
  p \\
\end{array}\right)\frac{1}{\alpha !}\left(
\begin{array}{c}
  p'+q' \\
  p' \\
\end{array}\right))^{\rho/\varrho})^{1/\rho} =$$
$$(\sum_{\beta \beta'}\frac{\vert
A_{\beta,\beta'}\vert^{\rho}}{\beta !(p!p'!)^{\rho
-1}})^{1/\rho}(\sum_{\gamma \gamma'}\frac{\vert B_{\gamma
\gamma'}\vert^{\rho}}{\gamma !(q!q'!)^{\rho
-1}})^{1/\rho}=\|A\|\|B\|$$ due to $\rho/\varrho =\rho-1$

To show the inequality $\lambda(p, p', q, q')\|A\|\|B\| \leq
\|A\bigodot B\|$ let us consider $$X=\{(A,B): A\in
Mat(p,p';F),\|A\|=1, B\in Mat(q,q';F),\|B\|=1\}$$, which is a
compact set
 in the corresponding finite dimensional vector space, and the  continuous map $(A,B)\mapsto A\bigodot B$.
 The image of $X$ with respect to this map is a compact set which doesn't contain zero vector
 due to Proposition 2. Let $\lambda(p, p', q, q')> 0$ stand for the distance between zero vector and
 this image set with respect to the corresponding $\rho$- norm. So $\lambda(p, p', q, q')\leq \|A\bigodot B\|$
 for any $(A,B)\in X$ and due to Proposition 2 one has
 $$\lambda(p, p', q, q')\|A\|\|B\| \leq \|A\bigodot B\|$$ for any $A\in Mat(p,p';F), B\in Mat(q,q';F)$

 With respect to the ordinary product of matrices a result similar to $\|A\bigodot B\|\leq
  \|A\|\|B\|$ is not valid. But one can have the following result.

  {\bf Proposition 7.} The following inequality
   $$\|A(p,q)B(q,q')\|\leq (p!)^{2-1/\rho}(q'!)^{2/\rho -1}\|A\|\|B\|_{\varrho}$$
   is true.

 {\bf Proof.} Indeed due to the H\"{o}lder inequality one has
 $$\|A(p,q)B(q,q')\|^{\rho}=\sum_{\alpha,\alpha'}\frac{1}{\alpha
!(p!q'!)^{\rho -1}}\vert\sum_{\beta}A_{\alpha,\beta}B_{\beta,\alpha'}\vert^{\rho}
\leq \sum_{\alpha,\alpha'}\frac{1}{\alpha
!(p!q'!)^{\rho -1}}\sum_{\beta}\vert A_{\alpha,\beta}\vert^{\rho}(\sum_{\gamma}\vert
B_{\gamma,\alpha'}\vert^{\varrho})^{\rho/\varrho}=$$
$$\sum_{\alpha,\beta}\frac{\vert A_{\alpha,\beta}\vert^{\rho}}{\alpha
!(p!q!)^{\rho -1}}(\sum_{\gamma,\alpha'}\frac{\vert
B_{\gamma,\alpha'}\vert^{\varrho}}{\gamma !(q!q'!)^{\varrho
-1}}\gamma !)^{\rho/\varrho}(q!)^{\rho}(q'!)^{2-\rho}\leq
\|A\|^{\rho}\|B\|_{\varrho}^ {\rho/\varrho}(q!)^{2\rho
-1}(q'!)^{2-\rho}$$ as far as $\gamma !\leq q!$.

In particular case the following estimation is also true.

 {\bf Proposition 8.} For any nonnegative integer numbers $m, k,
q'$ and $h\in Mat_{n,n}(0,1;F)$, \\ $A\in Mat_{n,n'}(m+k,q';F)$
the following inequality
 $$\|(\frac{h^{(m)}}{m!}\bigodot E_k)A\|\leq \left(\begin{array}{c}
  m+k \\
  k \\\end{array}\right)\|h\|_{\varrho}^m\|A\|$$ is valid.

{\bf Proof.} Indeed  $$\|(\frac{h^{(m)}}{m!}\bigodot
E_k)A\|^{\rho}=\sum_{\alpha,\alpha'}\frac{1}{\alpha !(k!q'!)^{\rho
-1}}\vert((\frac{h^{(m)}}{m!}\bigodot
E_k)A)_{\alpha,\alpha'}\vert^{\rho}=\sum_{\alpha,\alpha'}\frac{1}{\alpha
!(k!q'!)^{\rho -1}}\vert\sum_{\beta}(\frac{h^{(m)}}{m!}\bigodot
E_k)_{\alpha,\beta}A_{\beta,\alpha'}\vert^{\rho}= $$

$$\sum_{\alpha,\alpha'}\frac{1}{\alpha
!(k!q'!)^{\rho
-1}}\vert\sum_{\beta}\frac{h^{\beta-\alpha}}{(\beta-\alpha)!^{1/\varrho}}
\frac{A_{\beta,\alpha'}}{(\beta-\alpha)!^{1/\rho}}\vert^{\rho}$$
as far as
$$(\frac{h^{(m)}}{m!}\bigodot
E_k)_{\alpha,\beta}=\frac{h^{\beta-\alpha}}{(\beta-\alpha)!}$$

Due to the H\"{o}lder inequality
$$(\sum_{\beta}\vert\frac{h^{\beta-\alpha}}{(\beta-\alpha)!^{1/\varrho}}
\frac{A_{\beta,\alpha'}}{(\beta-\alpha)!^{1/\rho}}\vert)^{\rho}\leq
\sum_{\beta}\frac{\vert
A_{\beta,\alpha'}\vert^{\rho}}{(\beta-\alpha)!}
\sum_{\beta}\frac{\vert
h^{\varrho(\beta-\alpha)}\vert}{(\beta-\alpha)!}^{\rho/\varrho}
=\sum_{\beta}\frac{\vert
A_{\beta,\alpha'}\vert^{\rho}}{(\beta-\alpha)!}
(\frac{\|h\|_{\varrho}^{m\varrho}}{m!})^{\rho/\varrho}$$ Therefore
$$\|(\frac{h^{(m)}}{m!}\bigodot
E_k)A\|^{\rho}\leq \sum_{\beta,\alpha'}\frac{\vert
A_{\beta,\alpha'}\vert^{\rho}}{\beta! ((m+k)!q'!)^{\rho
-1}}\left(\begin{array}{c}
  m+k \\
  k \\\end{array}\right)^{\rho -1}\sum_{\alpha}\left(\begin{array}{c}
  \beta \\
  \alpha \\\end{array}\right) \|h\|_{\varrho}^{m\rho}=$$
  $$\sum_{\beta,\alpha'}\frac{\vert A_{\beta,\alpha'}\vert^{\rho}}{\beta!
((m+k)!q'!)^{\rho-1}}\left(\begin{array}{c}
  m+k \\
  k \\\end{array}\right)^{\rho}\|h\|_{\varrho}^{m\rho}= \|A\|^{\rho}\left(\begin{array}{c}
  m+k \\
  k \\\end{array}\right)^{\rho}\|h\|_{\varrho}^{m\rho}$$

 In future $Mat_{n,n'}(R)=Mat(R)$ stands for the set of all
block matrices $A=(A(p,p'))_{p,p'}$ with blocks $A(p,p')\in
M_{n,n'}(p,p';R)$ for all nonnegative integers $p$, $p'$. In
future it is assumed that $M(p,p';R)$ is a subset of $Mat(R)$ by
identifying each $A(p,p')\in M(p,p';R)$ as the element of $Mat(R)$
which's all blocks are zero, may be, except for $(p,p')$ block
which is $A(p,p')$.

For any $A, B\in Mat(R)$ we  define $A\bigodot B=C\in Mat(R)$,
where for all nonnegative integers $p$, $p'$
$$C(p,p')=\sum_{q,q'}A(q,q')\bigodot B(p-q,p'-q')$$

The above Propositions show that $(Mat(R); +, \bigodot)$ is an
integral domain, when $R$ itself is an integral domain. Its
identity element will be $1\in Mat(R)$ whose all blocks are zero
except for $(0,0)$ block which is $1$ -the identity element of
$R$.

In the case of $R=F$ we define $\rho$- norm $\|A\|=\|A\|_{\rho}$
of $A=(A(p,p'))_{p,p'}\in Mat(F)$, whenever it has meaning, in the
following form.

{\bf Definition 3.} $$\|A\|=\sum_{p,p'}\|A(p,p')\|$$

{\bf Theorem 1'.} If $A,B\in Mat(F)$ and $\lambda \in F$ then

a)$\|A\|=0$ if and only if $A=0$,

b)$\|\lambda A\| =|\lambda |\|A\|$,

c) $\|A+B\|\leq \|A\|+\|B\|$.

d) $\|A\bigodot B\|\leq \|A\|\|B\|$

{\bf Proof.} Here is a proof of d). Due to Theorem 1 one has
$$\|A\bigodot B\|=\sum_{p,p'}\|(A\bigodot
B)(p,p')\|=\sum_{p,p'}\|\sum_{q\leq p,q'\leq p'}A(q,q')\bigodot
B(p-q,p'-q')\| \leq $$
$$\sum_{p,p'}\sum_{q\leq p,q'\leq p'}\|A(q,q')\bigodot
B(p-q,p'-q')\|\leq \sum_{p,p'}\sum_{q\leq p,q'\leq
p'}\|A(q,q')\|\|B(p-q,p'-q')\|\leq$$
$$\sum_{q,q'}\|A(q,q')\|\sum_{p,p'}\|B(p,p')\|=\|A\|\|B\|$$

 In future the expression
$Exp(A)$, whenever it has meaning, stands for
$$E+\frac{1}{1!}A+\frac{1}{2!}A^{(2)}+\frac{1}{3!}A^{(3)}+ . .
.=\sum_{i=0}^{\infty}\frac{1}{i!}A^{(i)}$$, $R[x]$ is the ring of
polynomials in variables $x_1,x_2,...,x_n$  over $R$,
$x=(x_1,x_2,...,x_n)\in M_{n,n}(0,1;R[x])$.

Now one can easily derive the following result from Proposition 6.

{\bf Corollary 1.} If $B$ and $C$ are such matrices from
$M_{n',n''}(R)$ that each column of them has only finite number
nonzero elements then for any  $A=A(k,1)\in M_{n,n'}(k,1;R)$, the
following equality
$$Exp(A)B\bigodot Exp(A)C=Exp(A)(B\bigodot C)$$ is true.

If $n'\neq 0$ and
$$\varphi(x)=(\varphi_1(x),\varphi_2(x),...,\varphi_{n'}(x))=$$
$$x^{(0)}M_{\varphi}(0,1)+\frac{1}{1!}x^{(1)}M_{\varphi}(1,1)+
\frac{1}{2!}x^{(2)}M_{\varphi}(2,1)+...\in M_{n,n'}(0,1;R)$$ then
 one can screen it in the form
$$\varphi(x)=Exp{(x)}M_{\varphi} $$, where $M_{\varphi}\in
Mat(R)$ with blocks $M_{\varphi}(p,p')$ such that
$M_{\varphi}(p,p')=0$ whenever $p'\neq 1$ and only finite number
of blocks $M_{\varphi}(p,1)$ are not zero. We call $M_{\varphi}$
the matrix of the polynomial map $\varphi(x)$. Of course, if
$n'=0$ then
$$\varphi(x)=x^{(0)}M_{\varphi}(0,0)+\frac{1}{1!}x^{(1)}M_{\varphi}(1,0)+
\frac{1}{2!}x^{(2)}M_{\varphi}(2,0)+...\in M_{n,0}(0,0;R)$$

Now to understand the meaning of the product $\bigodot$ let us
assume that $n'=0$. Consider homogenous polynomials  $P=\sum
a_{\alpha}x^{\alpha}$ and $Q=\sum b_{\beta}x^{\beta}$ of degree
$m$ and $l$, respectively. It is not difficult to see that in this
case $$M_{PQ}=M_{P}\bigodot M_{Q}$$

 The Bombieri's norm of a
polynomial  $P(t)=\sum_{i=0}^m a_it^i$ is defined (in [1]) by
$[P]_2=(\sum_{i=0}^m \left(\begin{array}{c}
  m \\
  i \\
\end{array}\right)^{-1} a_i^2)^{\frac{1}{2}}$. Let us evaluate our $2$-norm of the
matrix $M_P$ of the corresponding homogeneous polynomial
$\sum_{i=0}^m a_it^is^{m-i}=\frac{(t,s)^{(m)}}{m!}M_P$, where
$M_P$ is the column matrix with entries $(a_ii!(m-i)!)$:
$$\|M_P\|_2=(\sum_{i=0}^m\frac{1}{i!(m-i)!m!}(a_ii!(m-i)!)^2)^{\frac{1}{2}}=[P]_2$$
that is in this case our $2$-norm and Bombieri's 2-norm are same.

The most remarkable feature of Bombieri's 2-norm states that for
any polynomials $P$, $Q$ the inequality
$$\left(\begin{array}{c}
  m+k \\
  k \\
\end{array}\right)^{1/2}[PQ]_2\geq [P]_2[Q]_2$$ is true, where
$m=deg P$, $k=deg Q$.

 With respect to the corresponding
matrices this inequality is nothing than
$$\left(\begin{array}{c}
  m+k \\
  k \\
\end{array}\right)^{-1/2}\|M_P\|_2\|M_Q\|_2\leq \|M_P \bigodot M_Q\|_2$$
Therefore we can say that in Theorem 1 we have a generalization of Bombieri's inequality.

 {\bf Theorem 2.} The following equality
$$Exp(Exp(x)M_{\varphi})=Exp(x)Exp(M_{\varphi}) $$ is valid.

{\bf Proof.} Indeed taking into account the above Propositions one
has
$$Exp(Exp(x)M_{\varphi})=\sum_{m=0}^{\infty}\frac{1}{m!}(M_{\varphi}(0,1)+\frac{1}{1!}x^{(1)}M_{\varphi}(1,1)+
\frac{1}{2!}x^{(2)}M_{\varphi}(2,1)+...)^{(m)}=\sum_{m=0}^{\infty}$$
$$\frac{1}{m!}\sum_{\alpha=(\alpha_0,...,\alpha_k,...),\vert \alpha\vert=m}\frac{m!}{\alpha_0!\alpha_1!...\alpha_k!...}
(M_{\varphi}(0,1))^{(\alpha_0)}\bigodot(\frac{1}{1!}x^{(1)}M_{\varphi}(1,1))^{(\alpha_1)}\bigodot...\bigodot
(\frac{1}{k!}x^{(k)}M_{\varphi}(k,1))^{(\alpha_k)}\bigodot...=$$

$$\sum_{i=0}^{\infty}\frac{x^{(i)}}{i!}\sum_{(\alpha_0,...,\alpha_k,...),1\alpha_1+...+k\alpha_k...=i}
\frac{M_{\varphi}(0,1)^{(\alpha_0)}}{\alpha_0!}\bigodot
\frac{M_{\varphi}(1,1)^{(\alpha_1)}}{\alpha_1!}\bigodot...\bigodot\frac{M_{\varphi}(k,1)^{(\alpha_k)}}{\alpha_k!}
\bigodot...$$

$$=\sum_{i=0,j=0}^{\infty}\frac{x^{(i)}}{i!}(Exp{(M_{\varphi}(0,1))}\bigodot
Exp{(M_{\varphi}(1,1))}\bigodot...\bigodot
Exp{(M_{\varphi}(k,1))}\bigodot...)(i,j)=$$
$$\sum_{i=0,j=0}^{\infty}\frac{x^{(i)}}{i!}
(Exp{(M_{\varphi})})(i,j)=Exp((x))Exp((M_{\varphi}))$$

 Consider
$\psi(y)=(\psi_1(y),\psi_2(y),...,\psi_n(y))=Exp{(y)}M_{\psi}$,
where $ M_{\psi}(i,1)\in M_{n'',n}(i,1;R)$ and

$$(\varphi \circ \psi)(y)=(\varphi_1(\psi(y)),
\varphi_2(\psi(y)),...,\varphi_{n'}(\psi(y)))=Exp{(y)}M_{\varphi
\circ \psi}$$

The following result is about the matrix representation of the
composition $\varphi \circ \psi$ of polynomial maps $\varphi$ and
$\psi$

{\bf Theorem 3.} The following equality $$M_{\varphi \circ
\psi}=Exp(M_{\psi})M_{\varphi}$$ is valid.

{\bf Proof.} $$(\varphi \circ \psi)(y)=Exp{(y)}M_{\varphi \circ
\psi}=\varphi(\psi(y))=Exp{(\psi(y))}M_{\varphi}=$$
 $$Exp{(Exp(y)M_{\psi})}M_{\varphi}=(Exp(y)Exp(M_{\psi}))M_{\varphi}=
Exp(y)(Exp(M_{\psi})M_{\varphi})$$ which implies that
$$M_{\varphi \circ
\psi}=Exp(M_{\psi})M_{\varphi}$$

{\bf Remark 2.} The equality $M_{\varphi \circ
\psi}=Exp(M_{\psi})M_{\varphi}$ indicates that  in our case the
real generalization of the ordinary product of matrices should be
the following binary operation $*$:
$$M_{\psi}*M_{\varphi}=Exp(M_{\psi})M_{\varphi}$$ as far as $M_{\psi}*M_{\varphi}$
 coincides with the ordinary product of matrices $M_{\psi}M_{\varphi}$ whenever $\psi$,
 $\varphi$ are linear maps.

 In a simple  case, when $\varphi(x)=
\frac{x^{(k)}}{k!}A(k,1)$,$\psi(x)= \frac{x^{(l)}}{l!}B(l,1)$ are
homogenous polynomial maps then due to Proposition 5 one has
$$\varphi(\psi(x))=\frac{(\psi(x))^{(k)}}{k!}A(k,1)=\frac{(\frac{x^{(l)}}{l!}B(l,1))^{(k)}}{k!}A(k,1)=
\frac{x^{(lk)}}{(lk)!}\frac{B(l,1)^{(k)}}{k!}A(k,1)$$ that is in
this case $$M_{\varphi \circ \psi}=\frac{B(l,1)^{(k)}}{k!}A(k,1)$$
Therefore  the following result is valid.

{\bf Corollary 2.} For any natural $m$ and $\varphi(x)=
\frac{x^{(k)}}{k!}A(k,1)$, where $A(k,1)=A\in Mat_{n,n}(k,1)$, the
following equality
$$M_{\varphi^{(m)}}=\frac{1}{k^{m-1}!}A^{(k^{m-1})}...\frac{1}{k^2!}A^{(k^2)}\frac{1}{k!}A^{(k)}A$$
is true, where
$\varphi^{(m)}(x)=\varphi(\varphi(...\varphi(x)...))$

The next result can be considered as a generalization of Theorem
1.

{\bf Theorem 4.} The following equality
$$Exp(Exp(M_{\psi})M_{\varphi})=Exp(M_{\psi})Exp(M_{\varphi})$$ is valid.

{\bf Proof.} Consider any polynomial map
$\xi(z)=(\xi_1(z),\xi_2(z),...,\xi_{n'}(z))=Exp{(z)}M_{\xi} \in
M_{n',n'''}(0,1;R)$. Due to $(\xi\circ\varphi)\circ\psi=
\xi\circ(\varphi\circ\psi)$ and Theorems 1,2 one has
$$M_{(\xi\circ\varphi)\circ\psi}=Exp(M_{\psi})M_{\xi\circ\varphi}=Exp(M_{\psi})(Exp(M_{\varphi})M_{
\xi})=(Exp(M_{\psi})Exp(M_{\varphi}))M_{ \xi},$$
$$M_{\xi\circ(\varphi\circ\psi)}=Exp(M_{\varphi\circ\psi})M_{\xi}=Exp(Exp(M_{\psi})M_{\varphi})M_{\xi}$$ and therefore
$$Exp(Exp(M_{\psi})M_{\varphi})=Exp(M_{\psi})Exp(M_{\varphi})$$

{\bf Corollary 3.} For any  $A=A(p,1)\in Mat_{n,n'}(p,1)$ and
$B=B(q,1)\in Mat_{n',n''}(q,1)$ the following equality
$$Exp(\frac{A^{(q)}}{q!}B)=Exp(\frac{A^{(q)}}{q!})Exp(B)$$
is true, in particular, if $A=A(1,1)\in Mat_{n,n}(1,1)$ is a
nonsingular matrix then $$(Exp(A))^{-1}=Exp(A^{-1})$$

{\bf Remark 3.} Usage of the introduced product makes power series
in many variables similar to the power series in one variable. It
allows to consider analogous problems for them considered in the
case of the power series in one variable and may be useful in
exploring analytical maps in many variables. For example, due to
Proposition 8 power series
$\sum_{m=0}^\infty\frac{x^{(m)}}{m!}A_m$ converges whenever
$\|x\|_{\varrho}<\frac{1}{r}$, where $A_m\in Mat_{n,n'}(m,q';F)$
and
$r=\overline{\lim}_{m\rightarrow\infty}\|A_m\|_{\rho}^{\frac{1}{m}}$.

\begin{center}{References}\end{center}

[1] O'Bryant, Kevin. "Bombieri Norm." From Math World--A Wolfram
Web Resource, created by Eric W. Weisstein.

http://mathworld.wolfram.com/BombieriNorm.html

 \end{document}